\documentclass[11pt,equ]{article}
\newtheorem{lemma}{Lemma}[section]
\newtheorem{prop}{Proposition}[section]
\newtheorem{theorem}{Theorem}[section]

\newcommand{\eqn}[1]{(\ref{#1})}

\newcommand{\hsp}{\hspace{\parindent}}

\newcommand{\bsq}{\vrule height .9ex width .8ex depth -.1ex}

\newcommand{\RR}{{\Bbb R}}
\newcommand{\QQ}{{\Bbb Q}}

\newcommand{\ZZ}{{\Bbb Z}}

\newcommand{\sgn}{{sgn}}

\newcommand{\beql}[1]{\begin{equation}\label{#1}}
\newcommand{\eeq}{\end{equation}}

\catcode`\@=11
\renewcommand{\section}{
        \setcounter{equation}{0}
        \@startsection {section}{1}{\z@}{-3.5ex plus -1ex minus
        -.2ex}{2.3ex plus .2ex}{\large\bf}
        }
\catcode`@=12
\makeatletter
\def\eqalignno#1{\displ@y \ta {\bf s} kip\@centering
  \halign to\displaywidth{\hfil$\@lign\displaystyle{##}$\ta {\bf s} kip\z@skip
    & $\@lign\displaystyle{{}##}$\hfil\ta {\bf s} kip\@centering
    & \llap{$\@lign##$}\ta {\bf s} kip\z@skip\crcr
    #1\crcr}}
\makeatother

\makeatletter
\def\@sect#1#2#3#4#5#6[#7]#8{\ifnum #2>\c@secnumdepth
     \def\@svsec{}\else 
     \refstepcounter{#1}\edef\@svsec{\csname the#1\endcsname.\hskip .75em }\fi
     \@tempskipa #5\relax
      \ifdim \@tempskipa>\z@ 
        \begingroup #6\relax
          \@hangfrom{\hskip #3\relax\@svsec}{\interlinepenalty \@M #8\par}%
        \endgroup
       \csname #1mark\endcsname{#7}\addcontentsline
         {toc}{#1}{\ifnum #2>\c@secnumdepth \else
                      \protect\numberline{\csname the#1\endcsname}\fi
                    #7}\else
        \def\@svsechd{#6\hskip #3\@svsec #8\csname #1mark\endcsname
                      {#7}\addcontentsline
                           {toc}{#1}{\ifnum #2>\c@secnumdepth \else
                             \protect\numberline{\csname the#1\endcsname}\fi
                       #7}}\fi
     \@xsect{#5}}
\def\@theorem#1#2{\it \trivlist \item[\hskip \labelsep{\bf #1\ #2.}]}
\input amssym.def
\input amssym.tex

\makeatother
\begin{document}
\begin{center}
{\large \bf Cyclic Systems of Simultaneous Congruences} \\ \bigskip

{\large {\em Jeffrey C. Lagarias}} \\
Department of Mathematics\\
University of Michigan \\
Ann Arbor, MI 48109-1109 \\
{\tt lagarias@umich.edu} \bigskip \\

(September 27, 2008) \\
\vspace*{1\baselineskip}
{\em Abstract}
\end{center}
This paper considers the cyclic system of $n \ge 2$ simultaneous congruences
$$
r  \left( \frac{\prod_{k=1}^n q_k}{q_i} \right)  
\equiv s~~( \bmod ~|q_i| ) ~, \quad 1 \leq i \leq n~,
$$
for fixed nonzero integers $(r,s)$ with $r>0$ and $(r,s)=1$. 
It shows there are only  finitely many solutions in positive integers 
$q_i \geq 2$, with $\gcd(q_1q_2 \cdots q_n, s)=1$ and obtains sharp bounds on the maximal size of solutions 
for almost all $(r,s)$. 
 The extremal solutions for $r=s=1$ are related to Sylvester's sequence
$2, 3, 7, 43, 1807,... .$ If the positivity condition on the
integers $q_i$ is dropped, then for $r=1$ 
these systems of congruences, taken $(\bmod ~|q_i|)$,
have infinitely many solutions, while for $r \ge 2$ they have finitely many solutions.
The problem is reduced to studying integer
solutions of the family of Diophantine equations
$$
r\left( \frac{1}{x_1} + \frac{1}{x_2} + \cdots + \frac{1}{x_n} \right) - \frac{s}{x_1 x_2 \cdots x_n} = m,
$$
depending on three parameters $(r, s, m)$. \\

\noindent

\setlength{\baselineskip}{1.0\baselineskip}
%
%
%
\section{Introduction}
\hsp
 Consider  the cyclic system of $n$ simultaneous congruences
\beql{100}
\left(\frac{q_1 q_2 \ldots q_n}{q_i} \right) \equiv 1 ~(\bmod ~q_i ) ~, 
\quad 1 \leq i \leq n ~,
\eeq
to be solved in integers $q_i \geq 2$. 
Sequences $(q_1, ..., q_n)$ with $2 \le q_1 \le  \cdots \le q_n$ satisfying \eqn{100} were named
 {\em Giuga sequences} by Borwein et al. \cite{BBBG96}, who
 related such sequences  to a conjecture of Giuga \cite{Gi50} on primality.
For $n = 2$ there are no solutions to
 \eqn{100}, but for $n =  3$ this
system has the unique solution
\begin{eqnarray}
2 \cdot 3 & \equiv & 1 \quad ( \bmod ~ 5) \nonumber \\
\label{102}
3 \cdot 5 & \equiv & 1 \quad ( \bmod ~ 2 )   \\
5 \cdot 2 & \equiv & 1 \quad ( \bmod~ 3) ~.  \nonumber
\end{eqnarray}
For $n \ge 3$ this equation has the solution
$(q_1, q_2, ..., q_n) = (u_1, u_2, ..., u_{n-1}, u_{n} -2)$
where $u_n$ are 
defined by the recursion $u_0=1$ and $u_{n+1}= (\prod_{i=1}^{n} u_i) +1$.
This sequence $u_n$  starts $1, 2, 3, 7, 43, 1807...$ and
grows doubly-exponentially in $n$. It is often called 
{\em Sylvester's sequence}, after 
work of J. J. Sylvester \cite{Syl80} in 1880. 
(However in  Knuth, Graham and Patashnik \cite{GKP94} this sequence is denoted
$e_n$ and its terms are called {\em Eulerian numbers}.)
 We show, 
 as a special case of Theorem~\ref{th11} below,
  that the solution above gives the maximal possible value of $q_n$
 in any Giuga sequence of length $n$.

In this paper we study solutions of the generalized system of
cyclic congruences
\beql{101}
r \left(\frac{q_1 q_2 \ldots q_n}{q_i} \right) \equiv s ~(\bmod ~q_i ) ~, 
\quad 1 \leq i \leq n ~,
\eeq
in which $r,s$ are nonzero  integers with $\gcd(r, s)=1$,  and 
we restrict to  solutions satisfying the greatest common divisor condition 
$
\gcd(q_1q_2 \cdots q_n, s)=1.
$
This  gcd condition is equivalent to the $q_k$ being pairwise 
relatively prime, as shown at the end of \S2.  
Without loss of generality  we reduce to the case $r >0$ by
multiplying the congruences by $-1$ if necessary; we allow $s$ to be positive
or negative. We also permit some variables $|q_i|=1$,
and call a solution  nontrivial if at least two $|q_i| \ge 2$.
We consider two situations: (1) the variables $q_j$ are restricted
to be positive integers; (2) the variables $q_i$ are nonzero integers.

The following result shows there are finitely many nontrivial positive integer 
solutions to  systems
of  simultaneous congruences \eqn{101} 
satisfying the gcd condition,  and gives a bound on their
size which is often sharp.

%
\begin{theorem}~\label{th11}
Let $r, s$ be nonzero integers with $r  >0$ and $\gcd(r, s)=1$. 
For each $n \geq 2$, there are only a finite number of  solutions  
in positive integers $(q_1, q_2, ..., q_n)\in (\ZZ_{+})^n$ 
to the cyclic system of $n$
simultaneous congruences
\beql{103}
r \left(\frac{\prod_{k=1}^n q_k }{q_i}\right) \equiv s~~(\bmod ~q_i ) 
~, \quad  1 \leq i \leq n. ~
\eeq
that satisfy the side conditions   (i) at least two $q_i \ge 2$,
and (ii)  
$\gcd(q_1q_2\cdots q_n, s)=1.$
Let the sequence $u_n(r)$ be given by $u_1(r)=r+1$ and
\beql{105}
u_{n+1}(r) = r\left( \prod_{i=1}^n u_i(r) \right) + 1. 
\eeq
Then each such solution to the cyclic system 
satisfies the upper bounds
\beql{104bb}
\max \{q_i \}  \leq  \left\{ \begin{array}{ll} 
\max \{ u_{n}(r)  - s -1,   s^2\} & ~\mbox{if}~~~ s>0,\\
\max \{ u_n(r)  -s -1, |s|\} & ~\mbox{if}~~~ s<0.
\end{array}
\right.
\eeq
For fixed $r,s$ the  upper bound  $u_n(r) - s-1$ is  attained for all
$n$ having   $u_{n}(r)  > s^2$. 
\end{theorem}

The case $r=s=1$ covers the case of Giuga sequences.
This theorem allows some moduli $q_i=1$ to occur in the
congruences; the corresponding congruence $(\bmod ~1)$
is  then automatically satisfied. Such moduli can be eliminated,
reducing the number of variables $n$ to cases where all $q_i \ge 2$,
retaining at least two such variables.
The two side conditions on solutions are necessary for finiteness, because  for  $n=2$,
$r=1$ and every  $s\ge 1$, there are an 
 infinite set of nontrivial positive solutions $\{(q_1, q_2)= (s, ks): k \ge 1\}$. 
 In Theorem~\ref{th11} 
 the case $s=1$ is excluded by  side condition (i)  that $q_1, q_2 \ge 2$, while
 all cases $s \ge 2$ are excluded by  side condition (ii) that $gcd(q_1q_2, s)=1$. 

The next two theorems concern  solutions to \eqn{101} allowing positive and
negative  integers. An
interesting feature here  is that for certain parameter values  there do exist
exist  infinitely many nontrivial solutions satisfying the side conditions (i),(ii).
As an example, for $n=3$ and $r=s=1$  the values
$(q_1, q_2, q_3) = (-k, k+1, k^2+k+1)$ for  $k \ge 2$
are  an infinite family of solutions. The following
result shows that whenever $r=1$ there are an infinite number of solutions.

%
\begin{theorem}~\label{th12}
Consider  the cyclic system of $n$
simultaneous congruences
\beql{111}
\frac{\prod_{k=1}^n q_k }{q_i} \equiv s~~(\bmod ~|q_i| ),
~, \quad  1 \leq i \leq n.
\eeq
where $s$ is nonzero.

(1) For each $n \ge 2$ this  system  has infinitely  integer solutions 
$(q_1, q_2, \cdots , q_n)\in (\ZZ \smallsetminus \{0\})^n$  with at
least two $|q_i| \ge 2$ and 
$$
\gcd(q_1 q_2 \cdots q_n, s) = 1.
$$
 
 (2) For each $n \ge2$ there  exists an integer $M_n^{\ast}$
 such that  when  $\gcd(s, M_n^{\ast})=1$,
  this  system   has infinitely  many integer solutions satisfying 
 $$
\gcd(q_1 q_2 \cdots q_n, s) = 1~~\mbox{and}~~~ \min\{ |q_i|\} \ge 2.
 $$ 
 An allowable value is  $M_n^{\ast} = u_1 u_2 \cdots u_n$,
 where $u_i= u_i(1)$ are terms in Sylvester's sequence. 
\end{theorem}

\noindent In the second part of this
result the  proof  determines the minimal values $M_2^{\ast}=1$ and $M_3^{\ast}=2$.
It does not determine the minimal value for $n \ge 3$, but it might be that 
that the general minimal value is $M_n :=gcd(2, n+1)$.

The next result shows that  in the  remaining cases $r \ge 2$ there are 
always a  finite number of integer
solutions, and obtains an  upper bound on their size.

%
\begin{theorem}~\label{th13}
Let  $r\ge 2$ and  $s$ be integers with $\gcd(r, s)=1$.
Then the  cyclic system of $n$ simultaneous congruences
\beql{104}
r \left(\frac{\prod_{k=1}^n q_k }{q_i}\right) \equiv s~~(\bmod ~|q_i| ) 
~, \quad  1 \leq i \leq n ~.
\eeq
 has only finitely  many integer solutions 
 $(q_1, q_2, \cdots, q_n) \in (\ZZ \smallsetminus \{0\})^n$ having \\
  $\gcd(q_1 q_2 \cdots q_n, s)=1$.
All such solutions satisfy the bound
\beql{104a}
\max\{ q_i \} \le \left( r(n+1)\right) ^{2^{n-1}}
+|s|.
\eeq
\end{theorem}

The upper bound  \eqn{104a} is far from tight; a slightly
better upper bound, more complicated to state,  is given in Theorem~\ref{th51}. 
It might be that the upper bound
of Theorem~\ref{th11} actually gives 
 the extremal bound for positive and negative variables, at least for
 $n$ large enough (depending on $r$ and$|s|$).

The proofs of Theorems \ref{th11}-\ref{th13} are based on 
reducing solutions  of the cyclic congruences to 
solutions of a family of 
Diophantine equations depending on 
three parameters $(r,s,m)$, namely
\beql{107a}
r\left( \frac{1}{x_1} + \frac{1}{x_2} + \cdots + \frac{1}{x_n} \right) - \frac{s}{x_1 x_2 \cdots x_n} = m,
\eeq
In  Lemma~\ref{le21} we show  that  a cyclic congruence  solution $(q_1, ..., q_n)$ satisfies \eqn{107a} 
for some integer $m$, a fact
noted by Borwein et al. \cite{BBBG96}. The condition $\gcd(q_1q_2 \cdots q_n, r|s|)=1$
leads to  a one-to-one correspondence of solutions
$(q_1, q_2, \cdots, q_n)$. In the other direction, all solutions
of \eqn{107a} (without the gcd restriction)
give solutions to the cyclic congruence system \eqn{101},
but some integer solutions $(q_1, q_2, \cdots q_n)$
to the cyclic congruence system \eqn{101} (without the gcd restriction) may not arise this way.

The main body of this  paper studies  integer solutions of
 the Diophantine equation \eqn{107a}, and does not 
impose any gcd restrictions on the variables. 
Many special cases of this equation have been previously considered in the literature,
and we discuss them below. 
To these results, this  paper supplies necessary and sufficient conditions when 
these equations have infinitely many integer solutions, given
in Theorem~\ref{th61}. 
We also establish  finiteness bounds  on the sizes
of the solutions that apply to all other cases 
of the Diophantine equation \eqn{107a}, given in  
Theorem ~\ref{th31} and Theorem~\ref{th51}.

The existence of an infinite number of integer solutions 
to the equation \eqn{107a} essentially  traces back to
the special case $m=0$, which we treat in \S4.
We use it to obtain a characterization of which parameter values
$(r,s,m)$  allow an infinite
number of solutions, which applies  more generally to the affine algebraic
hypersurface obtained from \eqn{107a} by clearing denominators.
This hypersurface is given by
\beql{120a}
r \left(\prod_{i=1}^n x_i \right) \left( \frac{1}{x_1} + \cdots + \frac{1}{x_n} \right)  + s = 
m \left( \prod_{i=1}^n x_i\right).
\eeq

The following theorem gives  necessary
and sufficient conditions on $(r,s, m)$ for this Diophantine equation to
have infinitely many integer solutions.

%
\begin{theorem}~\label{th14}
Let $r, s$  be  nonzero integers, with $r\ge 1$. Then for $n \ge 2$ the
affine algebraic hypersurface 
\beql{121}
r \left(\prod_{i=1}^n x_i \right) \left( \frac{1}{x_1} + \cdots + \frac{1}{x_n} \right)  + s = 
m \left( \prod_{i=1}^n x_i\right)
\eeq
defined over $\ZZ$ has infinitely many  integer solutions
$(x_1, x_2 , \cdots , x_n) \in (\ZZ \smallsetminus \{0\})^n$  if and only if $r= 1$ and one of 
the following conditions hold:  (i)   $|m| \le n-2$ and $s$ is arbitrary; (ii) $m=n-1$ and $s=1$; or (iii) $m=-(n-1)$ and $s= (-1)^{n-1}$.
\end{theorem}

For all parameter values $(r, s, m)$ the equation
\eqn{121} has infinitely many rational solutions $(x_1, \cdots , x_n) \in \QQ^n$.
Namely, if  we fix variables $x_2, \cdots, x_n$ to take rational values, then 
the remaining variable $x_1$ is
determined by a linear equation, so is rational. Thus  
the parameter restrictions of Theorem~\ref{th14}  are   a consequence of requiring
integrality of solutions.

Various authors have studied special cases of the Diophantine equation \eqn{107a},
often arising as a byproduct of studies on the Diophantine equation
\beql{108bb}
\frac{1}{x_1} + \cdots + \frac{1}{x_n} =\frac{a}{b}, 
\eeq
which encodes the problem of representing  $ \frac{a}{b}$ as a
sum of Egyptian fractions, where $0 < \frac{a}{b} \le 1.$ 
This work starts with  J. J. Sylvester \cite{Syl80} in 1880.
For $\frac{a}{b} =  1$ the double-exponential solution $u_n$  is
 related to the
problem of determining
\beql{109}
F_n : = 
\max \{ q_n : \sum_{i=1}^n \frac {1}{q_i} = 1 
\quad \mbox{with integers $q_i \geq 1$}\}~.
\eeq
It is well-known that the answer to \eqn{109} is
\beql{110}
F_n = u_n -1~,
\eeq
in which  $u_i = u_i (1)$ 
is the Sylvester sequence.
In  1921  Kellogg \cite{Kel21}  conjectured the equality \eqn{110}, 
which was  proved in 1922 independently  by Curtis \cite{Cur22} and
Takenouchi \cite{Tak22}. 
This bound was  reproved in 1950 by 
Erd\H{o}s \cite{Erd50}, in  the course of a more general investigation of
Egyptian fractions which raised new questions, cf.  Schinzel \cite{Sc02}.
  Another proof is given by  Soundararajan \cite{So05}. 
The extremal solutions to this problem turn out to have last variable $q_n$
expressible in terms of the preceding ones by
$q_n =q_1 q_2 \cdots q_{n-1}$, which yields  the Diophantine equation 
\beql{110b}
\left( \frac {1}{q_1} + \frac {1}{q_2} + ... + \frac {1}{q_{n-1} } \right)+
\frac {1}{q_1 q_2 ...q_{n-1}} = 1.
\eeq
This corresponds to the case $(r,s,m)= (1, -1, 1)$ in \eqn{107a}. 
The equation \eqn{110b} was directly studied by Brenton and Hill \cite{BH88} in
connection with complex surface singularities, and they determined
a complete list of positive solutions  for $n \le 8$.
Sylvester's sequence  also appears in connection with  
certain extremal lattice point problems, see
for example  Zaks, Perles and Wills \cite{ZPW82}
and Hensley \cite{Hen83}, who are concerned with
the maximal volume of a lattice simplex in $\RR^n$ containing
exactly $k$ lattice points. The author encountered
cyclic congruences in studying variants of such lattice point problems,
see Lagarias and Ziegler \cite{LZ91}. For recent work on a related lattice point
problem see Nill \cite{Ni07}.

An interesting open problem, that we do not consider, 
 is that of estimating the number of solutions 
to the  cyclic congruences \eqn{101} of size $n$, satisfying the side conditions (i), (ii)
for given parameters $(r, s)$.
This includes as special cases that of  counting
the number of Giuga sequences of length $n$, and of counting
the number of different Egyptian fractions of length $n$ that add up to $1$. 
Other  unsolved problems about Giuga sequences are listed 
in Borwein et al. \cite[Sect. 4]{BBBG96} and in  Borwein and Wong \cite{Bor97}.

The contents of this paper are as follows. In \S2 we give the relation of the 
cyclic congruence problem to solutions of the Diophantine equation \eqn{107a}.
In \S3 we consider positive solutions to
\eqn{107a} and use these to prove Theorem~\ref{th11}. Here a 
crucial ingredient (Proposition~\ref{pr31})
slightly extends the method of Erd\H{o}s \cite{Erd50}. In \S4 we study solutions to \eqn{107a}
having $m=0$, and determine when integer solutions occur in this case, with and
without the side condition $gcd(x_1 \cdots x_n, s)=1$, and use the results to
prove Theorem~\ref{th12}. In \S5 we show that when $m \ne 0$ there are only
finitely many solutions of \eqn{107a} having $\min\{ |x_i|\} \ge 2$, and give
a bound on their size, independent of $m$. We then use this result to prove
Theorem~\ref{th13}. In \S6 we determine for the case $r=1$ the parameters
$(s,m)$  for which the Diophantine equation \eqn{107a} has infinitely many solutions, and 
we use this to prove Theorem ~\ref{th14}.

\paragraph{Acknowledgments.} Part of this work was done during a visit 
 to the Mathematical Sciences Research Institute,
Berkeley, which is  supported in part by NSF, and part at AT\& T Labs-Research. 
The author thanks R. Girgensohn  for bringing a result of Erd\H{o}s to his
attention, via  Knuth,
Graham and Patashnik \cite[Exercise 4.59]{GKP94}, and thanks E. Croot for discussions
on the number of solutions. The author is indebted to the reviewer for many useful comments.
This work was  supported in part by NSF grant DMS-0500555.

%
%

\section{Associated Diophantine equation}

Solutions to cyclic congruences correspond to solutions of an
associated 
Diophantine equation of type \eqn{107a} with variable $m$, a fact observed
by Borwein et al. \cite{BBBG96} for Giuga sequences.

%
\begin{lemma}~\label{le21}
Let $r, s$ be nonzero integers with $r> 0$ and $\gcd(r,s)=1$.
Then for each $n \ge 2$ the following conditions are equivalent.

(1)  The nonzero integers $(q_1, ..., q_n)$ with $\gcd(q_1q_2 \cdots q_n, s)=1$
satisfy 
the cyclic system of simultaneous congruences
\beql{201}
r \left(\frac{\prod_{k=1}^n q_k }{q_i}\right) \equiv s~~(\bmod ~|q_i| ) 
~, \quad  1 \leq i \leq n ~.
\eeq

(2) The nonzero integers  $(q_1, ..., q_n)$ with
$gcd(q_1q_2 \cdots q_n, s)=1$
satisfy the Diophantine equation
\beql{202}
r\left( \frac{1}{q_1} + \cdots \frac{1}{q_n}\right) - \frac{s}{q_1q_2 \cdots q_n} = m
\eeq
for some integer $m$.
\end{lemma}

\paragraph{Proof.} Suppose (1) holds, and write
\beql{203}
M := r \left( \prod_{k=1}^n q_k \right)
\left( \frac{1}{q_1} + \cdots  +\frac{1}{q_n}\right) - s.
\eeq
Each $q_i$ divides the integer 
$M$, for it divides $\frac{r}{q_j}(\prod_{k=1}^n q_k)$ for each $j \ne i$, and
it divides $r \left(\frac{\prod_{k=1}^n q_k}{q_i}\right) - s$ by \eqn{201}. We conclude that
the least common multiple $[q_1, q_2, \cdots, q_n]$ of the $q_i$ divides $M$.

The cyclic congruence \eqn{201} and the fact $\gcd(q_1q_2 \cdots q_n, s)=1$ implies that
$r\left(\frac{q_1q_2 \cdots q_n}{q_i}\right)$ is invertible $(\bmod~|q_i|)$, whence
$(q_j, q_i)=1$ for all $j \ne i$. Thus the $q_i$ are pairwise relatively prime, so
$[q_1, q_2, \cdots, q_n]= q_1q_2 \cdots q_m$. Thus we can write
$$
M =m (q_1 q_2 \cdots q_n)
$$
for some integer $m$. Dividing \eqn{203}
by $q_1q_2 \cdots q_m$ yields the Diophantine equation \eqn{202}, which gives  (2). 

Suppose (2) holds, without imposing the gcd condition \eqn{101a}.
 Multiply \eqn{202} by $q_1q_2 \cdots q_n$ to obtain
$$
r \left( \sum_{i=1}^n \frac{q_1q_2 \cdots+ q_n}{q_i} \right) - s = mq_1q_2 \cdots q_n,
$$
Reducing this equation $(\bmod |q_i|) $ yields a solution to the cyclic congruence
\eqn{201} for $q_i$, which is (1), without imposing the gcd condition. Now imposing the
gcd condition gives (1).  $~~~\bsq$

\paragraph{Remarks.}
(1) This proof  shows  that 
solutions to the Diophantine equation \eqn{202} not satisfying  the gcd condition 
$gcd(q_1q_2\cdots q_n, s)=1$ 
still  give solutions to the corresponding cyclic congruence \eqn{201}, 
not satisfying  the gcd condition. 
However the converse is not true, for one may  take $n=2$
and $(r,s)=(1, -20)$, and then 
$(q_1, q_2)=(5, 25)$ has $\gcd(q_1q_2, s)= 5$  and satisfies the cyclic congruence
system \eqn{201} but not  the Diophantine equation \eqn{202}.\\

(2) All solutions to the cyclic congruences \eqn{201} with $(q_1\cdots q_n,  s) =1$
necessarily have
$$\gcd(mq_1 q_2 \cdots q_n, r) = 1.$$ 
To see this, consider  the associated Diophantine equation \eqn{202} and multiply by
all the $q_i$ to obtain
$$
r \left( \prod_{k=1}^n q_k \right)
\left( \frac{1}{q_1} + \cdots  +\frac{1}{q_n}\right) - s= m (\prod_{i=1}^n q_i).
$$
Reducing this equation $(\bmod~r) $ yields
$$
m\prod_{i=1}^n q_i \equiv s~(\bmod ~r).
$$
Since $s$ is invertible $(\bmod~m)$ we obtain $gcd(mq_1 q_2 \cdots q_n, r) = 1.$\\

(3) For the cyclic congurence \eqn{201}, the side condition 
$\gcd(q_1 q_2 \cdots q_n, s) = 1$ holds if and only if  all $q_i$ are pairwise
relatively prime, i.e. $\gcd(q_i, q_j)=1$ if $i \ne j$. 
The "only if" direction was shown in the proof of Lemma~\ref{le21}.
For the "if" direction, we prove the contrapositive. If 
$\gcd(q_1q_2\cdots q_n, s)>1$, there exists some prime $p|s$ 
with $p|  \gcd(q_i, s)$ for some index $i$.  Then 
the $i$-th cyclic congruence
$r \left(\frac{\prod_{k=1}^n q_k }{q_i}\right) \equiv s~~(\bmod ~|q_i| ) $
yields
$$
r \left(\prod_{k \ne i} q_k \right)\equiv 0~(\bmod~p),
$$
whence $\gcd (r, s)=1$ yields  $p|q_k$ for some $k \ne i$, so 
 $p| \gcd(q_i, q_k)$, and the $q_i$ are not pairwise relatively prime.

%
%

\section{Positive integer solutions: Proof of Theorem~\ref{th11}}

We treat the case of positive solutions $(q_1, ..., q_n)$ to the
cyclic congruence \eqn{101}, and prove Theorem~\ref{th11}.
 Lemma~\ref{le21}  reduces  this to questions about positive
integer solutions  of  the Diophantine equation
\beql{301}
r\left( \frac{1}{x_1} + \frac{1}{x_2} + \cdots  + \frac{1}{x_n} \right) - 
\frac{s}{ x_1 x_2 \cdots x_n}= z.
\eeq
The following result 
gives  bounds on the size of positive
solutions to this equation, without imposing any
gcd conditions on the solutions. 

%
\begin{theorem}~\label{th31}
Suppose $r, s$ are nonzero integers with $r \ge 1$ and $(r, s) =1$.
Then for  fixed $n \ge 2$ the
Diophantine equation \eqn{301}
 has only finitely many integer solutions 
$(x_1, x_2, \cdots, x_n, z) =(q_1, q_2, ..., q_n, m)$  satisfying
$1 \le q_1 \le q_2 \le \cdots \le q_n$, and the side conditions
\beql{301b}
  r q_1\cdots q_{n-1} > s, ~~~\mbox{if}~~s >0,
\eeq
\beql{301c}
   q_n > |s| ~~~\mbox{if}~~s <0,
\eeq
and allowing any  $m$.
Furthermore, when these conditions hold: 

(1) All solutions have $m \ge 1$ and satisfy the bounds 
\beql{302}
\max \{q_i \}  \leq  \left\{ \begin{array}{ll} 
u_{n}(r)  - 2  & ~\mbox{if}~~~ s>0,\\
u_n( r) +|s| -1  & ~\mbox{if}~~~ s<0.
\end{array}
\right.
\eeq
Here $u_n(r)$ is determined by the recursion
$u_1(r)= r+1$ and
\beql{303}
u_{n+1}(r) = r\left( \prod_{i=1}^n u_i(r) \right) + 1. 
\eeq

(2) If   in addition $q_n > s^2$, then there holds 
\beql{302b}
\max\{ q_i\} \le u_n(r)   -s -1.
\eeq
Furthermore, whenever $u_n(r) > s^2$,
 there  exist solutions for which $q_n= u_n(r) - s -1$, so the bound \eqn{302b} is 
 tight for all   $(r,s, n)$ for which  
$u_n(r) > s^2$. 
\end{theorem}

The condition $r q_1\cdots q_{n-1}>s$ is needed in the theorem when $s>0$
because 
the case  $r q_1\cdots q_{n-1}=s$ 
 can sometimes have  infinitely many positive solutions, when no gcd conditions
 are imposed. 
 These occur
when $\sum_{i=1}^{n-1} \frac{1}{q_i} = \frac{m}{r}$ for some positive integer $m$.
The condition  $q_n > |s|$ when $s<0$  is needed for a similar reason.

We will use the following bound for a Kellogg-type Diophantine
equation, which is an extension to all $r$ of  a result
of Erd\H{o}s \cite[T\'{e}tel 5]{Erd50} for 
the case $r=1$. (A sketch of Erd\H{o}s's  proof appears in
Knuth, Graham and Patashnik \cite[Exercise 4.59; cf.  p. 523]{GKP94}.)
%
\begin{prop}~\label{pr31}
Let $r \ge 1$ be an integer, and 
let $p_1, p_2, \cdots p_k$ be positive integers with
$k \ge 1$ such that
\beql{324}
\sum_{i=1}^k \frac{1}{p_i} + \frac{1}{\alpha_{k+1}} = \frac{1}{r}.
\eeq
in which $\alpha_{k+1}$ is a positive rational number
such that
\beql{325} 
\alpha_{k+1} \ge  \max \{p_1, p_2, \cdots , p_k \}. 
\eeq
Let $u_1( r)=r+ 1$ and $u_{k+1}(r) = r u_1(r) \cdots  u_k(r) + 1$.
Then $\alpha_{k+1}$ satisfies the bounds
\beql{326}
\alpha_{k+1} \le u_{k+1}( r) -1,
\eeq
\beql{327}
p_1p_2 \cdots p_k ( \alpha_{k+1}+1) \le u_1(r) u_2(r)\cdots u_k(r) u_{k+1}(r).
\eeq
\end{prop}

\paragraph{Proof.}
We prove the result by induction on $k \ge 1$. For the base case  $k=1$,
 \eqn{324} requires $p_1 \ge r+1$,
which  then implies 
$$
\alpha_2 \le r(r+1) = u_2(r) -1,
$$
since $\frac{1}{r+1} + \frac{1}{r(r+1)} = \frac{1}{r}$. 
The condition  $\alpha_2 \ge p_1$ requires $r+1 \le p_1 \le 2r$.
To complete the base case we must verify that 
$$
p_1 (\alpha_2 + 1) \le(r+1)(r^2+r+1) =  u_1(r) u_2(r).
$$
To see this, we solve \eqn{324} for $\alpha_2$,
obtaining
$\alpha_2= \frac{p_1 r}{p_1 -r}.$
Viewing this as a function of $p_1$ we set
$$
f(p_1):= p_1(\alpha_2 + 1) = p_1(\frac{p_1 r}{p_1 -r} + 1).
$$
We calculate  for $p_1>r$, 
$$
\frac{d^2f}{dp_{1}^2}  = \frac{2p_1r}{p_1-r} + \frac{2(p_1)^2 r}{(p_1-r)^3} \ge 0.
$$
Thus $f(p_1)$ is convex downwards on $r+1 \le p_1 \le 2r$ so its maximum occurs
at one of the endpoints of this interval. Now
$ f(2r) = 2r(2r+1) \le f(r+1) = (r+1)(r^2+r+1),$ for $r\ge1$, giving the result.

For the induction step, we assume it is true for all values $k-1$ or smaller and treat
the given $k \ge 2$, which concerns
\beql{327aa}
\frac{1}{p_1} + \cdots + \frac{1}{p_k} + \frac{1}{\alpha_{k+1}} = \frac{1}{r}.
\eeq
Now there  exists an integer $d \ge 1$ such that
$$
\frac{1}{\alpha_{k+1}}= \frac{1}{r} - \left( \sum_{i=1}^k \frac{1}{p_i} \right) = \frac{d}{r p_1p_2 \cdots p_k},
$$
hence
\beql{327a}
 \alpha_{k+1} = \frac{r p_1p_2 \cdots p_k}{d} \le  r p_1p_2 \cdots p_k.
\eeq
The induction step is treated in three cases.

\paragraph{Case 1.} $p_k > p_{k-1}$
and $\frac{1}{p_1} + \cdots + \frac{1}{p_{k-1}} + \frac{1}{p_k -1} \ge \frac{1}{r}$.\\

In this case there exists a rational $\beta_k$ with $p_k > \beta_k \ge p_k - 1 \ge p_{k-1}$ such that
$$
\frac{1}{p_1} + \cdots + \frac{1}{p_{k-1}}+ \frac{1}{\beta_k} = \frac{1}{r},
$$
The induction hypothesis for $k-1$ is satisfied for this equation, so we conclude,
$$
  \beta_k +1 \le u_k(r),
$$
and, since $p_k \le \beta_k+1$,  
\beql{328}
p_1 p_2 \cdots p_k \le p_1 p_2 \cdots p_{k-1} ( \beta_k +1) \le u_1(r) u_2(r) \cdots u_k(r).
\eeq
Combined with \eqn{327a} this  yields
$$
 \alpha_{k+1} \le r p_1 p_2 \cdots p_k \le r u_1(r) u_2(r) \cdots u_k(r) = u_{k+1}(r) -1.
$$
Consequently, using \eqn{328}, 
$$
p_1 p_2 \cdots p_k ( \alpha_{k+1} + 1) \le u_1(r) u_2(r) \cdots u_k(r) u_{k+1}(r),
$$
completing the induction step in Case 1.

\paragraph{Case 2.} $p_k = p_{k-1}$
and $\frac{1}{p_1} + \cdots + \frac{1}{p_{k-1}} + \frac{1}{p_k -1} \ge \frac{1}{r}$.\\

We may suppose $p_k \ge 4$, for the only case with $p_k \le 3$  has  
 $p_k=p_{k-1}=3$ with $k=2, r=1$ and $\alpha_{k+1}= 3$, which
 satisfies the theorem.
When $p_k \ge 4$ we have $\frac{1}{p_k -2} \le \frac{2}{p_k}= \frac{1}{p_k}+ \frac{1}{p_{k-1}}$,
so there exists $\beta_k >0$ such that 
\beql{328b}
\frac{1}{p_1} + \cdots + \frac{1}{p_{k-2}} + \left(\frac{1}{p_k -2} + \frac{1}{\beta_k}\right) = \frac{1}{r}.
\eeq
We claim that
\beql{328c}
(p_k)^2 \le (p_k-2)(\beta_k +1).
\eeq
Assuming this claim is proved, it implies $\beta_k > p_k$, hence the equation \eqn{328b}
satisfies the induction hypothesis for $k-1$. This  yields
$$
p_1 p_2 \cdots p_{k-2} p_{k-1}p_k \le p_1 p_2 \cdots p_{k-2}(p_k-2)(\beta_{k}+1) 
 \le u_1(r) u_2(r) \cdots u_{k}(r).
$$
Then using \eqn{327a} 
we obtain 
$$
 \alpha_{k+1} \le r p_1 p_2 \cdots p_k \le r u_1(r) u_2(r) \cdots u_k(r) = u_{k+1}(r) -1.
$$
Consequently
$$
p_1 p_2 \cdots p_k (  \alpha_{k+1} + 1) \le u_1(r) u_2(r) \cdots u_k(r) u_{k+1}(r),
$$
completing the induction step in case 2.

It remains to prove the claim \eqn{328c}. Subtracting  \eqn{327aa} from \eqn{328b}
and using $p_k=p_{k-1}$  gives
$$
\frac{1}{\beta_k} =\frac{1}{\alpha_{k+1}}+ \left(  \frac{2}{p_k} - \frac{1}{p_k - 2} \right).
$$
The Case 2 hypothesis gives
$$
\frac{1}{\alpha_{k+1}} \le \frac{1}{p_k -1} - \frac{1}{p_k} = \frac{1}{p_k(p_k -1)},
$$
so that 
$$
\frac{1}{\beta_k} \le 
\frac{1}{p_k (p_k-1)} +  \frac{2}{p_k} - \frac{1}{p_k - 2} =
\frac{(p_k)^2 - 4p_k +2}{p_k(p_k-1)(p_k-2)}.
$$
Setting $y=p_{k}$, this yields 
$$
\beta_k \ge \frac{y(y-1)(y-2)}{y^2-4y+2}.
$$
from which we obtain
$$
(p_k -2)(\beta_{k} +1) \ge \frac{y(y-1)(y-2)^2}{y^2-4y+2}.
$$
Now \eqn{328c} follows from the inequality
$$
 \frac{y(y-1)(y-2)^2}{y^2-4y+2} \ge y^2  ~~~\mbox{for}~~ y \ge 4,
$$
which is easily verified by clearing the denominator and simplifying.

\paragraph{Case 3.} $\frac{1}{p_1} + \cdots + \frac{1}{p_{k-1}} + \frac{1}{p_k -1} < \frac{1}{r}$.\\

For reasons as in case 2, we may suppose $p_k \ge 5$. 
In this case we  replace $p_{k}$ by $\tilde{p}_k := p_k-1$, and
 $\alpha_{k+1}$ by $\tilde{\alpha}_{k+1}$ satisfying
\beql{330}
\frac{1}{p_k} + \frac{1}{\alpha_{k+1}} = \frac{1}{p_k-1}+ \frac{1}{\tilde{\alpha}_{k+1}},
\eeq
The case 3 inequality guarantees that $\tilde{\alpha}_{k+1}$ is positive, which
ensures that $\tilde{\alpha}_{k+1} > \alpha_{k+1} \ge p_k$.
We claim that in addition
\beql{331}
 p_k( \alpha_{k+1} +1) \le (p_k-1) (\tilde{\alpha}_{k+1}+1).
\eeq
Assuming this claim is proved, we have obtained a new equation of the same size $k$, 
$$
\frac{1}{p_1} + \cdots + \frac{1}{p_{k-1}} + \left(\frac{1}{p_k -1} + \frac{1}{\tilde{\alpha}_{k+1}} \right) 
= \frac{1}{r}, 
$$
which has new denominators $(p_1, ..., p_{k-2}, p_{k-1}, \tilde{p}_k)$, which are the same or
smaller than the original system,
while $\tilde{\alpha}_{k+1}$ has increased. It  then suffices to prove the induction step for the
new equation because it would give
$$
\alpha_{k+1} \le \tilde{\alpha}_{k+1} \le u_{k+1}(r) -1
$$
$$ 
p_1  \cdots p_{k-1} p_k(\alpha_{k+1}+1) \le p_1 \cdots p_{k-1} (p_{k} - 1)(\tilde{\alpha}_{k+1} +1)
\le u_1(r) u_2(r) \cdots u_k(r) u_{k+1}(r),
$$
thus proving the induction step for the original equation. 
If the new system falls in case 1 or 2 we are done, while if it falls in case 3 we can repeat 
the reduction. Eventual
termination into case 1 or 2 must occur, because the sum of the denominators
decreases at each step. Thus case 3 will terminate.

It remains to prove the claim \eqn{331}.  Using the definition \eqn{330} of 
$\tilde{\alpha}_{k+1}$, we  express it in terms of $y=p_k$ and $\alpha_{k+1}$,
obtaining
$$
\tilde{\alpha}_{k+1}= \frac{y(y-1) \alpha_{k+1}}{y(y-1)- \alpha_{k+1}},
$$
noting that $y(y-1) > \alpha_{k+1}$ under the Case 3 hypothesis.
Substituting this in the inequality  \eqn{331} 
and clearing a positive denominator shows that  it is equivalent to 
$$
(\frac{y}{y-1}) \alpha_{k+1}^2 - (y+1)\alpha_{k+1} + y(y-1) \ge 0,
$$
when $\alpha_{k+1} \ge y \ge 5$. The left side is a quadratic polynomial in $\alpha_{k+1}$
having discriminant $D= (y+1)^2 - 4y^2$, and $D<0$  when $y >1$, hence
it is then positive for all 
real  $\alpha_{k+1}$
and the inequality  \eqn{331} follows. This completes Case 3. $~~~\bsq$\\

We use Proposition \ref{pr31} to establish Theorem~\ref{th31}.
%
\paragraph{Proof of Theorem~\ref{th31}.}

We first eliminate  all values $q_i=1$, which reduces the equation \eqn{301} to
an equation in fewer $x$-variables having the same form, with $z$
shifted by an integer, and without affecting the side conditions. 
The bound to be proved is nondecreasing  in $n$ and independent of $z$, so it suffices
to prove the upper bound for the new system, which has 
solutions satisfying   $2 \le q_1 \le  q_2 \le \cdots \le  q_n$.
We denote the 
associated  integer choice of the $z$ variable by $m$. Now we can rewrite a solution to \eqn{301} as
\beql{340}
\sum_{i=1}^n \frac{1}{p_i} + \frac{1}{\alpha_{n+1} }= \frac{m}{r},
\eeq
where we have 
$p_i :=  q_i$ and $\alpha_{n+1} := -\frac{r}{s} q_1 q_2 \cdots q_n,$
where $\alpha_{n+1}$ may be positive or negative, depending on the sign of $s$.

\paragraph{Case 1.} $s>0$. \\

 We claim that 
$r q_1 \cdots q_{n-1} > s$ 
 implies $m \ge 1$.
To show this,  note that we always have, for some 
 integer $b >0$, 
 \beql{340a}
\left(\sum_{i=1}^{n-1} \frac{1}{q_i}\right) + \frac{1}{q_n}\left( 1- \frac{s}{rq_1 q_2 \cdots q_{n-1}}
\right)  = \frac{b}{q_1\cdots q_{n-1} }+ \frac{1}{q_n}\left( 1- \frac{s}{rq_1 q_2 \cdots q_{n-1}}
\right) = \frac{m}{r}.
\eeq
The condition $rq_1 q_2 \cdots q_{n-1} >s$ implies that the coefficient of $\frac{1}{q_{n}}$ is
positive, so we may 
legitimately  solve this equation for $q_n$, obtaining
\beql{341}
q_n = \frac{ 1-  \frac{s}{rq_1 q_2 \cdots q_{n-1}}}{\frac{m}{r} - \frac{b}{q_1q_2 \cdots q_{n-1}}}
=\frac{rq_1q_2 \cdots q_{n-1} -s} {mq_1 \cdots q_{n-1}-rb}.
\eeq
Here  the  denominator of the fraction cannot vanish, since
$q_n < \infty$.
 If now $m \le 0$,
then the denominator of the fraction on the right in \eqn{341} would be
negative, while the numerator is positive, contradicting $q_n >0$. Thus $m \ge 1$,
and the claim is proved. 

Now  we  divide \eqn{340a} by $m\ge 1$ to obtain
\beql{342}
\left(\sum_{i=1}^{n-1} \frac{1}{mq_i}\right) + \left( \frac{1}{mq_n} -
 \frac{s}{mrq_1q_2 \cdots q_n}\right) = \frac{1}{r}.
\eeq
We view this as an instance of Proposition~\ref{pr31},
with  $k=n-1$,  setting   $p_i= mq_i$ for $1 \le i \le n-1$,
and 
$$
\frac{1}{\alpha_n} := \frac{1}{mq_n} -  \frac{s}{mrq_1q_2 \cdots q_n}
= \frac{1}{mq_n} (1 - \frac{s}{rq_1q_2 \cdots q_{n-1}}).
$$
Now we have 
\beql{343}
\alpha_n= m q_n \left(\frac{rq_1 q_2 \cdots q_{n-1}}{rq_1q_2 \cdots q_{n-1}-s}\right), 
\eeq
and  $rq_1 \cdots q_{n-1} >s$   yields 
\beql{344}
\alpha_n > mq_n \ge m q_{n-1}=p_{n-1} >0,
\eeq
so the hypotheses of Proposition~\ref{pr31} are satisfied for \eqn{342}.
The proposition then gives
\beql{345}
\alpha_n \le u_{n}(r)-1,
\eeq
\beql{346} 
q_1 q_2 \cdots q_{n-1}(\alpha_n+1) \le  u_1(r) u_2(r) \cdots u_{n}(r).
\eeq
Solving \eqn{343} for $q_n$  yields 
$$
q_n =  \frac{1}{m} \left(1- \frac{s}{rq_1q_2 \cdots q_{n-1}}\right) \alpha_n
\le  \frac{1}{m} \left(1- \frac{1}{rq_1q_2 \cdots q_{n-1}}\right)
(u_{n}(r)-1).
$$
Now \eqn{346} implies
\begin{eqnarray*}
\frac{s}{rq_1 \cdots q_{n-1} }&=& \frac{s (\alpha_n+1)}{rq_1 \cdots q_{n-1} (\alpha_n +1)} \\
&\ge& \frac{s(\alpha_n +1)}{ru_1(r) \cdots u_{n-1}(r) u_n(r)}
= \frac{s (\alpha_n +1)}{(u_{n}(r) -1) u_n(r)}.
\end{eqnarray*}
From \eqn{343} and \eqn{344} we then deduce
\begin{eqnarray}\label{347a}
q_n  &\le & u_n(r) - 1 - \frac{s (\alpha_n +1)}{ u_n(r)}  \nonumber \\
& \le & u_n(r) - 1 - \frac{s(m q_n)}{u_n(r)} \nonumber  \\
& \le & u_n(r) - 1 - \lceil \frac{s q_n}{u_n(r) } \rceil.
\end{eqnarray}
From this inequality we immediately deduce that for $q_n >1$, 
\beql{348}
q_n \le u_n(r) - 2.
\eeq
 We also deduce that, for $q_n > s^2$, there holds 
\beql{349}
q_n \le u_n(r) - s -1.
\eeq
Suppose otherwise, so that $q_n \ge u_n(r) - s.$ Then \eqn{347a} gives,
$$
u_n(r) -s \le q_n \le u_n(r) - 1 - \lceil s- \frac{s^2}{u_n(r)}\rceil = u_n(r) -s -1 +
 \lfloor \frac{s^2}{u_n(r)}\rfloor = u_n(r) -s -1,
$$
a contradiction which establishes \eqn{349}.

\paragraph{Case 2.} $s <0.$\\
Now the left side of \eqn{301} has every term positive, which implies
that $z=m >0$, so again $m \ge 1$. Dividing by $m$ we obtain
\beql{351}
\left( \sum_{i=1}^{n-1} \frac{1}{mq_i} \right) + \frac{1}{mq_n} + \frac{|s|}{mr q_1 \cdots q_n} = \frac{1}{r}.
\eeq
By hypothesis $q_n > |s|$, and we consider  a new system that replaces $q_n$ by $q_n -|s|$.
Now there is an 
integer $b$ such that
\beql{352}
\left( \sum_{i=1}^{n-1} \frac{1}{mq_i} \right) + \frac{1}{m(q_n-|s|)} + 
\frac{b}{mr q_1 \cdots q_{n-1}(q_n- |s|)} = \frac{1}{r}.
\eeq
Subtracting \eqn{351} from this equation and rearranging terms yields
\begin{eqnarray*}
\frac{b}{mr q_1 \cdots q_{n-1}(q_n- |s|)} &=& \frac{1}{mq_n} +
 \frac{|s|}{mr q_1 \cdots q_n} - \frac{1}{m(q_n-|s|)} \\
&=& 
\frac{rq_1 \cdots q_{n-1}(q_n - |s|) + |s|(q_n-|s|) -r q_1 \cdots q_n}{mrq_1 \cdots q_{n-1}q_n(q_n-|s|)} \\
&=& 
\frac{|s|(q_n - |s| - r q_1 \cdots q_{n-1})}{mrq_1 \cdots q_{n-1}q_n(q_n-|s|)}
\end{eqnarray*}
Comparing both sides yields
\beql{353}
b= |s| \left( 1 - \frac{rq_1 q_2 \cdots q_{n-1} + |s|}{q_n}\right).
\eeq
We claim $b \le 0$, which is the same as
\beql{353a}
q_n \le rq_1q_2\cdots q_{n-1} + |s|.
\eeq
To show this, note that there is a positive integer $b'$ such that
\beql{353b}
\frac{1}{r} - \left( \sum_{i=1}^{n-1} \frac{1}{mq_i} \right) = \frac{b'}{mrq_1 \cdots q_{n-1} } \ge 
\frac{1}{mr q_1 \cdots q_{n-1}}.
\eeq
By \eqn{351} the  left side of this expression equals
$$
\frac{1}{mq_n}(1+ \frac{|s|}{r q_1 \cdots q_{n-1}}) = \left(\frac{r q_1 \cdots q_{n-1} + |s|}{mq_n}\right)
\frac{1}{mrq_1 \cdots q_{n-1}}.
$$
Comparison with \eqn{353b} yields \eqn{353a}, proving the claim.

We  treat two subcases, $b=0$ and $b <0$.

\paragraph{Subcase 2.1.} $b=0$.\\

In this subcase we have
$$
\left( \sum_{i=1}^{n-1} \frac{1}{mq_i}\right) + \frac{1}{m(q_n - |s|)}= \frac{1}{r}.
$$
This is a system of the form of Proposition \ref{pr31} for $k=n-1$, after permuting the
terms to take $\alpha_n = \max\{ mq_i, m(q_n-|s|)\}$. We conclude from the
proposition that
$$
q_n-|s| \le m(q_n - |s|) \le \alpha_n \le u_n(r) - 1.
$$
Then we deduce
$$
q_n \le (q_n- |s|) +|s| \le u_n(r) + |s| -1.
$$

\paragraph{Subcase 2.2.} {\em  $b <0$.}\\

In this subcase we have
$$
\left( \sum_{i=1}^{n-1} \frac{1}{mq_i}\right) + \frac{1}{m(q_n - |s|)} -
\frac{|b|}{mr q_1q_2\cdots q_{n-1}(q_n - |s|)}= \frac{1}{r},
$$
which is a system of the  form \eqn{301}
with new variables $(q_1^{'}, \cdots , q_n^{'}): = (q_1,\cdots, q_{n-1}, q_n -|s|)$
and $z=m$, 
and with new parameters  $(r', s') = (r, |b|)$.
We claim this system falls under Case 1,  taking  
$s' = |b| >0$,
possibly after permuting the variables. 
We must verify that two side conditions hold, namely $q_n^{'}= q_n -|s| >0$ and 
\beql{357}
rq_1^{'} q_2^{'} \cdots q_{n-1}^{'} > |b|.
\eeq
The first of these holds  by hypothesis.
For the second, we observe that \eqn{357} is equivalent to 
\beql{358}
rq_1 q_2 \cdots q_{n-1} +|s| > |b| + |s|.
\eeq
But now  \eqn{353} gives (since both $b, s <0$),
$$
|b| + |s|  = \frac{|s|(rq_1 q_2 \cdots q_{n-1} +|s|)}{q_n} = 
(rq_1 q_2 \cdots q_{n-1} +|s|) \frac{|s|}{q_n} < rq_1 q_2 \cdots q_{n-1} +|s|,
$$
which verifies  \eqn{358}.
Thus the  Case 1 hypotheses are met for the new system.

We now  apply the Case 1 inequality in the form \eqn{348} to obtain 
$$
q_n - |s| \le \max\{ q_i^{'}\} \le u_n(r) - 2,
$$
and this  yields
$$
q_n \le u_n(r) + |s| -2,
$$
in Subcase 2.2.

Combining the two subcases,  we conclude that  in Case 2 one always has 
$$
q_n \le u_n(r) + |s| -1 = u_n(r) -s-1,
$$
as required.

To finish the proof it remains to verify  the tightness of the upper bounds 
$u_n(r)-s-1$ for given  nonzero $r,s$ 
with $r>0$ and $gcd(r,s)=1$, for those $n$ having solutions with  $q_n> s^2$
if $s>0$ and $q_n>|s|$ if $s<0.$ 
We show the existence of solutions
with $q_n =u_n -s -1.$
One verifies, by induction on 
$n$,  that
$$
(x_1, x_2, \cdots, x_n) = (u_1(r), u_2(r), \cdots, u_{n-1}(r), u_n(r) -s+1)$$
 gives
a solution to \eqn{301} with $z=1$  for all pairs $(r, s)$ whenever 
$u_n(r) - s +1 \ge 2$. The key property is that 
$$
\frac{1}{r} - \left( \sum_{i=1}^{n-1} \frac{1}{u_i(r)} \right) = \frac{1}{ru_1(r) u_2(r) \cdots u_n(r)},
$$
and
$$
- \frac{1}{ru_1(r) u_2(r) \cdots u_{n-1}(r)} +\frac{ 1} {ru_1(r) \cdots u_{n-1}(r) - s} =
 \frac{s}{ru_1(r) \cdots u_{n-1}(r) ( u_n(r) -s - 1)}.
$$
 Note also that $gcd(q_1 q_2 \cdots q_n, r) =1$ always holds  for these solutions,
using $gcd(r, s)=1$, but $\gcd(q_1 q_2 \cdots q_n, s) >1$ may occur  for certain $|s|>1$.
$~~~\bsq$\\

%

We now deduce Theorem ~\ref{th11} from Theorem ~\ref{th31}.

\paragraph{Proof of Theorem~\ref{th11}.}
To bound the size of  a positive solution $(q_1, ..., q_n)$ to the
cyclic congruence, we first drop all variables $q_i=1$, which reduces
to a cyclic congruence in fewer variables, but necessarily at least
two variables, since at least two $q_i \ge 2$. The desired upper
bounds \eqn{104bb} are all nondecreasing functions of $n$ so it suffices
to treat the smaller problem. Thus we may assume all $q_i \ge 2$,
and we may  reorder the variables so that
 $2 \le q_1 \le  q_2 \le \cdots \le  q_n$, with $n \ge 2$.

We then use Lemma~\ref{le21} to convert the cyclic congruence solution
$(q_1, ..., q_n)$ to a solution
of a Diophantine equation \eqn{301} with some value of $z$,
and we may presume this solution is large enough to satisfy the
side conditions \eqn{301b} and \eqn{301c} in Theorem~\ref{th31}. The equation we
consider is then
\beql{304}
\left( \sum_{i=1}^{n} \frac{1}{q_i} \right)  - \frac{s}{rq_1q_2 \cdots q_n} = \frac{m}{r}.
\eeq
By the discussion after Lemma~\ref{le21} the hypotheses 
$\gcd(q_1 \cdots q_n, s) =1$ implies $\gcd(q_i, q_j)=1$ if $i \ne j$
and $\gcd(m q_1 \cdots q_n, r)=1$. Thus we must have 
$2 \le q_2 < q_3 < \cdots < q_n $. 

We show finiteness of the number of solutions. 
Consider first the case $s<0$. 
Theorem~\ref{th31} establishes  finiteness of 
solutions having   $q_{n}>s$.  Finiteness of solutions 
having $q_n \le s$ is immediate, since they then satisfy
$2 \le q_1< q_2 < \cdots < q_n \le s$. 
Now consider the remaining cases  $s>0$.  
Theorem~\ref{th31} also establishes finiteness of the number of
solutions, for fixed $s>0$, that satisfy   $rq_1 \cdots q_{n-1} > s.$
For the remaining cases with $r q_1 \cdots q_{n-1}\le s$, we establish
finiteness of the admissible solutions (those having $\gcd(q_1 \cdots q_n, s)=1$) 
by obtaining an upper bound for $q_n$, namely $q_n <  s.$
We first show there are no admisible
solutions in the equality case   $rq_1 q_2 \cdots q_{n-1} =s$. For  $|s|\ge 2$ this holds 
since $(q_1\cdots q_n, s)=1$ and $(r, s)=1$ by hypothesis, and for $s= \pm 1$
it holds because $q_1 \ge 2$. 
(The equality case is a critical case, for a single solution 
to it would yield infinitely many positive solutions to \eqn{103}, since $q_n$ is  unconstrained.)
 It remains  to  treat  cases where $rq_1 \cdots q_n < s$.
Now we have, by  \eqn{341}, that the Diophantine equation can be solved
for $q_n$, 
\beql{361}
q_n :=\frac{rq_1q_2 \cdots q_{n-1} -s} {mq_1 \cdots q_{n-1}-rb},
\eeq
where
$$
\sum_{i=1}^{n-1} \frac{1}{mq_i} = \frac{b}{mrq_1 \cdots q_{n-1}}.
$$
The numerator in \eqn{361} does not vanish, and the finiteness of
$q_n$ requires the denominator be nonzero, whence
$$
q_n \le  |rq_1q_2 \cdots q_{n-1} -s| = s - r q_1 \cdots q_n < s.
$$
Thus finiteness follows. 

The explicit bounds follow from Theorem ~\ref{th31} .
The condition $u_{n-1}(r) > s$ implies that
$r u_1(r) \cdots u_{n-1}(r) > s$  and that 
$$
u_n(r) = (u_{n-1}(r )-1) u_{n-1}(r) \ge s(s+1) > s^2,
$$
so the desired conclusion $q_n \le u_{n}(r) - s +1$
follows in these cases. $~~~\bsq$\\

%
%
%

\section{General integer solutions with  $m=0$ : Proof of Theorem~\ref{th12}}
\hsp

In the remainder of the paper we study general integer solutions $(x_1 , \cdots, x_n)$ to the
Diophantine equation 
\beql{400}
r\left( \frac{1}{x_1} + \frac{1}{x_2} + \cdots + \frac{1}{x_n} \right) - \frac{s}{x_1 x_2 \cdots x_n} = m,
\eeq
in which $r,s $ are integers with $r >0$ and $\gcd(r,s)=1$ . In this section we treat the special
case $m=0$. We show that infinitely many solutions occur when $r=1$, and then 
characterize when an infinite number of solutions exist satisfying extra gcd  conditions.
We apply  these results  to prove  Theorem~\ref{th12}. 

%

\begin{theorem}~\label{th41}
Suppose $r, s$ are nonzero integers with $r>0$ and $\gcd(r,s)=1$. Consider the Diophantine equation
\beql{401}
r\left( \sum_{i=1}^n \frac{1}{x_i} \right) = \frac{s}{x_1x_2 \cdots x_n}
\eeq
Then:

(1) If $r > 1$ this equation has   no integer solutions $(x_1, \cdots , x_n)$.

(2) If $r=1$ this equation has infinitely many integer solutions.
Moreover there
are infinitely many integer solutions having $\min \{ |x_i| \} \ge 2$.
\end{theorem}

%
\paragraph{Proof.}
(1) Suppose $r \ge 2$ and that an integer solution $(x_1, \cdots , x_n)$ exists. Then
$$
\sum_{i=1}^n \frac{1}{x_i} = \frac{s}{r x_1 \cdots x_n}.
$$
However
$$
\sum_{i=1}^n \frac{1}{x_i} = \frac{b}{x_1x_2 \cdots x_n}= \frac{br}{x_1x_2 \cdots x_n},
$$
for some integer $b$. Thus $s= br$, which contradicts $\gcd(r, s) =1$ since
$r \ge 2$.

(2) Now let $r=1$. Define the sequence $\{v_k(m) : k \ge 0\}$ for each $m \ge 1$ by
the recurrence $v_1(m)=m$ and 
\beql{411}
v_k(m) := \left(v_1(m) v_2(m) \cdots v_{k-1}(m)\right) +1,
\eeq
which gives the sequence $(m, m+1, m(m+1) + 1, \cdots).$
One proves, by induction on $n \ge 2$, that
\beql{412}
-\frac{1}{v_1(m)} + \left( \sum_{i=2}^{n-1} \frac{1}{v_i(m)}\right) = \frac{-1}{v_1(m) v_2(m) \cdots v_{n-1}(m)}.
\eeq
Now we obtain
\begin{eqnarray}\label{413}
-\frac{1}{v_1(m)} + \left( \sum_{i=2}^{n-1} \frac{1}{v_i(m)} \right) + \frac{1}{v_n(m)+ (s-1)}
&=& - \frac{1}{v_1(m) \cdots v_{n-1}(m)} + \frac{v_1(m) \cdots v_{n-1}(m) +s}  \nonumber\\
&=& \frac{-s}{v_1(m) \cdots v_{n-1}(m) ( v_n(m)-s)}.
\end{eqnarray}
It follows that, for $m > |s|+1 \ge 2$, 
$$
(x_1, \cdots , x_{n-1}, x_n) :=  (-v_1(m), v_2(m), \cdots , v_{n-1}(m), v(n)(m) -s)
$$
satisfies
$$
\sum_{i=1}^n \frac{1}{x_i} = \frac{s}{x_1 \cdots x_n}
$$
since $x_1= x_1(m)<0$ and all other $x_i=x_i (m)>0$, and 
 $\min\{ |x_i|\} \ge 2.$  $~~~\bsq$\\

The next result shows that if we require that infinitely many integer 
solutions satisfying   the extra condition
$gcd(x_1x_2 \cdots x_n, s)=1$ required in the cyclic congruence,
then the set of parameters allowing such solutions
narrows  slightly. By Theorem~\ref{th41} we need
only consider the case that $r=1$.
 In part (1) of  the following result we include  integer solutions
in which some variables  $x_i= \pm 1$.

%

\begin{theorem}~\label{th42}
Let $s$ be a nonzero integer, and consider the Diophantine equation 
\beql{421}
\frac{1}{x_1} + \frac{1}{x_2} + \cdots +  \frac{1}{x_n} = \frac{s}{x_1x_2 \cdots x_n},
\eeq
where $n \ge 2$.

(1) This equation has infinitely many integer solutions $(x_1, ..., x_n)$ with at
least two $|x_i| \ge 2$ and 
 $\gcd(x_1x_2 \cdots x_n, s)=1$ for those $(n,s)$ such that  $\gcd(s, M_n)=1$, where
 $M_{2k}=1$ and $M_{2k+1}=2$. For all remaining $(n,s)$, 
 where  $n$ is odd and $s$ is even, this equation has no integer solutions.

(2) For each $n \ge 2$ there is a finite modulus $M_n^{\ast}$ such that
whenever  $\gcd(s, M_n^{\ast}) =1$ this equation 
 has infinitely many integer solutions satisfying $\min\{ |x_i| \} \ge 2$
 with $\gcd(x_1x_2, \cdots x_n, s)=1$.
One can take $M_n^{\ast} = u_1 u_2 \cdots u_{n},$  where the $u_i$ belong to
Sylvester's sequence. 
\end{theorem}

%
\paragraph{Proof.}

(1) {\em Necessity.} We show there are
no solutions with $\gcd(x_1 x_2 \cdots x_n, s)=1$ when  $n \equiv 1~(\bmod~2)$ and $s$ is even.
The gcd condition implies that all $x_i$ are odd, and multiplying \eqn{421} by
$x_1x_2 \cdots x_n$ yields
\beql{431}
\sum_{i=1}^n \frac{x_1x_2 \cdots x_n}{x_1} = s.
\eeq
Now the left side is odd, being a sum of $n$ odd terms, and the right side is even,
a contradiction. 

{\em Sufficiency.} Suppose $n \ge 5.$ Choosing $x_n=1, x_{n-1}=1$ leads to 
an equation of the same form with $n-2$ variables, and with $s$ replaced by $-s$.
In this way we reduce to the cases $n=2$ with any $s$ and to $n=3$ with an
odd $s$, and it remains to show infinitely many solutions having 
$\gcd(x_1 \cdots x_n, s) =1$ exist in these cases. For $n=2$ we have the infinite
family of solutions
$(x_1, x_2) = (-m, m+s)$ giving
$$
-\frac{1}{m} + \frac{1}{m+s} = \frac{s}{(-m)(m+s)},
$$
taking $m \ge |s|+2$. It now suffices to restrict to the
arithmetic progression  $m \equiv 1~(\bmod ~|s|)$. For $n=3$ we consider the family
$(x_1, x_2, x_3) =(-m, m+1, m(m+1) + s)$, for $m \ge |s|+2$. Here $s$ is odd, so if we again
choose $m \equiv 1~(\bmod~|s|)$, then $(-m(m+1)(m+s), s) = 1$. 

(2) Theorem~\ref{th41} exhibited for each $n$ and each fixed  $s \ne0$ an infinite family of solutions
to \eqn{431} having $\min \{ |x_i(m) |\} \ge 2$. Now the sequence $m \equiv 1 ~(\bmod~s)$
will have the required property $gcd(x_(m) \cdots x_n(m), s) = 1$, as long as
$(x_1(1)\cdots x_n(1), s) =1$, since the congruence 
$x_j(1+ks) \equiv x_j(1)~(\bmod ~|s|) $ is easy to establish. 
Now we have $x_j(1) \equiv v_j(1) =u_j ~(\bmod ~s)$, whence
$$
x_1(1)\cdots x_n(1)= (v_1(1) v_2(1) \cdots v_{n-1}(1)(  v_n(1) -s) 
\equiv u_1 u_2 \cdots u_n (\bmod ~s).
$$
Thus if we take $M_{n}^{\ast} = u_1 u_2 \cdots u_n$, then
$\gcd(M_n^{\ast}, s)=1$ will imply $(q_1(1+jm)\cdots q_n(1+jm), s)=1$, for all $j$.
$~~~\bsq$\\

\paragraph{Remark.}  
One might  conjecture that the  minimal allowable value
of $M_n^{\ast}$ in Theorem ~\ref{th42}  equals $ M_n= gcd (n+1,2)$. 
The necessity part of the proof above showed that 
one can take $M_2^{\ast}=M_2=1$ and $M_3^{\ast} =M_3 = 2$, confirming
this conjecture in these cases.

%
\paragraph{Proof of Theorem~\ref{th12}.}
(1)  For  $r=1$, the cyclic congruence \eqn{103} becomes \eqn{111}.
 We show for each $n$ and all nonzero $s$ the latter
cyclic congruence   has infinitely many nontrivial 
solutions satisfying the gcd condition.  Theorem~\ref{th42} (1), 
together with Lemma~\ref{le21}, shows 
on choosing $m=0$ that there are infinitely many nontrivial solutions
to \eqn{103}
when $r=1$ and $\gcd(s, M_n)=1$, having 
$\min\{ |q_i|\} \ge 2$.   This handles all $s$ when $n$ is even.
If $n$ is odd, we choose one variable $q_n=1$. Eliminating this
variable, whose cyclic congruence $(\bmod ~q_n)$
is trivially satisfied, yields  a cyclic congruence in
$(n-1)$ variables with the same $(r, s)=(1,s)$, which now has
infinitely many nontrivial solutions satisfying
$\gcd(q_1\cdots q_{n-1}, s)=1$ by Theorem~\ref{th42}(1) since $n-1$ is even.\\

(2) This follows similarly from Theorem~\ref{th42} (2) together with Lemma~\ref{le21}.
$~~~\bsq$

%
%
%

\section{General integer solutions: Proof of Theorem~\ref{th13}}
\hsp

In this section we prove a finiteness theorem 
on the number of integer 
solutions to  \eqn{107a}  with $m \ne 0$, with  all variables $|x_i| \ge 2$, 
but with no gcd condition on solutions, and apply this result to 
prove Theorem~\ref{th13}.

%
\begin{theorem}~\label{th51}
Let $r, s$ be integers with $r>0$ and $\gcd(r,s)=1$.
Suppose $m \ne 0$ is fixed. Then the Diophantine equation
\beql{501}
r\left( \frac{1}{x_1} + \frac{1}{x_2} + \cdots  + \frac{1}{x_n} \right) - 
\frac{s}{  x_1 x_2 \cdots x_n}= m
\eeq
has finitely many  integer solutions  satisfying
$$
\min\{ |x_i|\} \ge 2.
$$
All such solutions satisfy the bound
\beql{503}
\max\{ |x_i| \} \le r^{2^{n-1}}\left(\prod_{j=1}^{n-1} (n+2-j)^{2^{n-1-j}}\right) +|s|.
\eeq
\end{theorem}

%
\paragraph{Proof.}
We first note that the hypotheses imply $\gcd(mx_1 \cdots x_n, r)=1$.
This follows by multiplying by $x_1 \cdots x_n$ and then reducing $(\bmod~r)$,
obtaining
$$
m x_1 x_2 \cdots x_n \equiv -s ~(\bmod~r),
$$
and the gcd result follows since $\gcd(r,s)=1$. 

Without loss of generality reorder the variables
so that $2 \le |x_1| \le |x_2 \le \cdots \le |x_n|$. We rewrite \eqn{501} as
\beql{504}
\frac{1}{x_1} + \cdots + \frac{1}{x_n} - \frac{s}{rx_1 \cdots x_n} = \frac{m}{r}.
\eeq
Now set
\beql{505}
R_i := \left( \frac{1}{x_i} + \frac{1}{x_{i+1}} +\cdots + \frac{1}{x_n} \right) - \frac{s}{rx_1 \cdots x_n}
\eeq
Now \eqn{504} gives
$$
|R_1| = \frac{|m|}{r} > 0.
$$
Assuming that $|x_n| >|s|$ we have
$$
|R_1| \le \frac{1}{|x_1|}\left( n + \frac{|s|}{r|x_2 \cdots x_n|} \right) \le \frac{1}{|x_1|}(n+1).
$$
This implies
\beql{506}
|x_1| \le \frac{1}{R_1}(n+1) \le \frac{r}{|m|}(n+1) \le r(n+1).
\eeq
Now \eqn{504} yields
\beql{507}
R_i = \frac{m}{r} -\left( \frac{1}{x_1} + \cdots + \frac{1}{x_{i-1}} \right)= \frac{m_i}{r x_1 x_2 \cdots x_{i-1}}
\eeq
with 
\beql{508}
m_i := m x_1 x_2 \cdots x_{i-1} - r x_1 x_2 \cdots x_{i-1}\left( \sum_{j=1}^{i-1} \frac{1}{x_j}\right).
\eeq 
We claim that $m_i \ne 0$. This follows since
$$
m_i \equiv mx_1 x_2 \cdots x_{i-1}~(\bmod~r),
$$
and we have shown $(m x_1 \cdots x_n, r) =1$. Thus we obtain
\beql{509}
|R_i| \ge \frac{1}{r|x_1 \cdots x_{i-1}|}.
\eeq
Combining this with the definition \eqn{505} yields
$$
\frac{1}{r|x_1 \cdots x_{i-1}|}  \le |R_i| \le \frac{1}{|x_i|}\left(n-i+1 + \frac{|s x_i|}{r |x_1 \cdots x_n|}\right)
$$
and this gives, for $2 \le i \le n-1$
\beql{511}
|x_i| \le r|x_1 \cdots x_{i-1}|(n-i+2) \le
\eeq
To bound  the last variable we use 
$$
R_n = \frac{m_n}{rx_1 \cdots x_{n-1}}= \frac{1}{x_n} \left( 1- \frac{s}{rx_1 \cdots x_{n-1}} \right) 
$$
which can be solved for $x_n$ to give
$$
x_n = \frac{rx_1 \cdots x_{n-1} -s} {m_n}.
$$
Since $m_n \ne 0$, this yields
\beql{512}
|x_n| \le r|x_1 x_2 \cdots x_{n-1} |+ |s|.
\eeq
Now \eqn{511} yields, by induction on $i \ge 1$, for $1 \le i \le n-1$, 
$$
|x_i| \le r ^{2^{i-1}} \left(\prod_{j=1}^{i-1} (n+2-j)^{2^{i-1-j}}\right) (n+2-i)
$$
with the base case given by \eqn{506}. Finally, \eqn{512} now gives
$$
|x_n| \le r ^{2^{n-1}}\left( \prod_{j=1}^{n-1}  (n+2-j)^{2^{n-1-j}} \right)+|s|,
$$
completing the proof. $~~~~\bsq$.\\

We now deduce Theorem~\ref{th13} from this result.

%
\paragraph{Proof of Theorem~\ref{th13}.}
We have $r \ge 2$. By Lemma ~\ref{le21} this corresponds to a
Diophantine equation of the type \eqn{501} with some value of $m$.
We eliminate variables $x_i= \pm 1$,  reducing
to a similar system with smaller $n$, and with a new $s' = \pm s$,
and a new value $m'$.  If the reduced system has 
$n=0$, then  all $|x_i|=1$. If  the reduced system has $n=1$,
i.e. all but one variable have $|x_i|=1$, say for $i \ge 2$, we
obtain the  reduced system
$$
\frac{r}{x_1}- \frac{s}{x_1}=m.
$$
Since $r \ge 2$ and $\gcd(r,s)=1$ we have $r-s \ne 0$, and the only choices of
$x_1$ giving  integer $m$ have $|x_1| $ dividing $ |r-s|$, so $|x_1| \le r + |s|,$
and \eqn{104a} holds.
Finally, if the reduced system has $n \ge 2$ then this system cannot have $m'=0$, because
Theorem ~\ref{th41}(1) says there are no solutions to the
resulting Diophantine equation with $m=0$. Thus $m \ne 0$, and 
the bound of Theorem~\ref{th51} applies. This bound \eqn{503} is stronger than
what is needed, since 
$$
r ^{2^{n-1}}\left( \prod_{j=1}^{n-1}  (n+2-j)^{2^{n-1-j}} \right)+|s|
\le \left( r(n+1)\right)^{2^{n-1}} +|s|,
$$
which gives \eqn{104a}.  $~~~~\bsq$\\

%
%
%

\section{General integer solutions: Proof of Theorem~\ref{th14}} 
\hsp

We now determine for arbitrary $m$  and $\gcd(r, s)=1$ when the Diophantine
equation \eqn{501} has infinitely many integer solutions, and apply this result to prove 
Theorem~\ref{th14}.  
The following result classfies
when an infinite number of solutions exists for $r=1$. 

%

\begin{theorem}~\label{th61}
Suppose that $r= 1$, and $s$ is nonzero. Then the Diophantine equation
\beql{601}
\left( \frac{1}{x_1} + \frac{1}{x_2} + \cdots + \frac{1}{x_n} \right) - \frac{s}{x_1 x_2 \cdots x_n} =m
\eeq
has infinitely many  integer solutions  if and only if  
(i)   $|m| \le n-2$ and $s$ is arbitrary, or (ii) $m=n-1$ and $s=1$, or (iii) $m=-(n-1)$ and $s=(-1)^{n-1}$.
\end{theorem}

\paragraph{Proof.}
{\em Necessity.}  We may suppose $1 \le |x_1| \le |x_2| \le \cdots \le |x_n|$. 
We first show there are finitely many
solutions when $|m|\ge n$. Set 
\beql{601b}
R := \frac{1}{x_1} + \frac{1}{x_2} + \cdots + \frac{1}{x_n} 
\eeq
so that \eqn{601} becomes 
\beql{602}
R+ \frac{s}{x_1x_2 \cdots x_n} = m.
\eeq
If  there were an infinite number of solutions to \eqn{602}, at least one has some $|x_i| \ge 6|s|$. 
For this solution $|\frac{s}{x_1x_2 \cdots x_n} | \le \frac{1}{3} $, so that 
$$
| R| \le \sum_{i=1}^n \frac{1}{|x_i|} \le (n -1) + \frac{1}{6|s|} \le  n- \frac{5}{6}.
$$
However \eqn{602} gives
$$
|R| \ge |m| - |\frac{s}{x_1x_2 \cdots x_n}| \ge |m| - \frac{1}{6}.
$$
Combining these yields   $|m| \le n - \frac{2}{3} < n,$ so $|m| \le n-1.$

 It remains to treat the cases $m = \pm(n-1)$.
 If  $|m|=(n-1)$, and 
 $|x_{n-2}| > 1$, then $|R| < n-2 + \frac{1}{2} + \frac{1}{6 |s|} \le n - \frac{4}{3},$
 while \eqn{602} gives $|R| \ge (n-1) - \frac{|s|}{6 |s|} = n - \frac{7}{6}$, a contradiction.
Thus we must have all $|x_i|=1$ for $1 \le i \le n-2$. For  the case $m=n-1$, we must have
all these $x_i=1$ and \eqn{601} simplifies to 
$$
 \frac{1}{x_n} - \frac{s}{x_n} =0.
$$
This has infinitely many solutions if and only if $s=1$. For the case  $m= -(n-1)$,
all $x_i=-1$ for $1 \le i \le n-1$ and \eqn{601} simplifies to 
$$
\frac{1}{x_n} - (-1)^{n-1} \frac{s}{x_n} =0.
$$
This equation has infinitely many integer solutions if and only if $s= (-1)^{n-1}$. \\

{\em Sufficiency}. We choose $x_n= \pm 1$ so that $\sgn (x_n)= \sgn (m)$. In that
case the equation becomes
\beql{631}
\left( \frac{1}{x_1} + \frac{1}{x_2} + \cdots + \frac{1}{x_{n-1} }\right) - 
\sgn(m) \frac{s}{x_1 x_2 \cdots x_{n-1}}
 = \sgn(m) (|m|-1).
\eeq
By multiplying this equation by $\sgn(x_n)$ we then get an equation of the
same form with one fewer variable, with $s' = \sgn(m) s$, and with  a right side decreased by
one in absolute value. 

Continuing in this way, since $|m| \le n-1$, we eventually arrive at a system with
$n' = n-|m|$ variables and right hand side value $m'=0$. 
For $n' = 1$ the system is of the form
$$
\frac{1}{x_1} - \frac{s'}{x_1} = 0,
$$
(where $s'= \pm s$) which has infinitely many solutions if and only if $s'=  1.$ This
corresponds to the two cases above.

For $n'=2$, it has the form
$$
\frac{1}{x_1} + \frac{1}{x_2} - \frac{s'}{x_1x_2} = 0
$$
These have the infinite family of solutions, for $m > |s'|+2$, 
$(x,_1, x_2) = (m, - (m-s')) $ if  $s'>0$, and
$(x_1, x_2) = (m, -(m+s'))$ if  $s'<0$. 
For $n'=3$ it has the form
$$
\frac{1}{x_1} + \frac{1}{x_2} + \frac{1}{x_3} - \frac{s'}{x_1x_2 x_3} =0.
$$
This has the family of solutions, for $m > |s'|+2$ 
$(x_1, x_2, x_3)= (m, -(m+1), -m(m+1) -s') $ 
if $s'>0$, 
$(x_1, x_2, x_3)= (-m,  m+1,  m(m-1)+s')$ 
if  $s'<0$. 

The cases with $m=0$ and $n \ge 4$  
 we can set two variables $x_n=1, x_{n-1}=-1$
and reduce to an equation with two fewer variables and still with $m'=0$.
Continuing this way, we reduce to an equation \eqn{301} with $n=2$ or $n=3$ variables,
 having $m=0$, which has infinitely many solutions by the
constructions above. $~~\bsq$\\

We now deduce Theorem~\ref{th14}.

\paragraph{Proof of Theorem~\ref{th14}.}
By hypothesis $r \ge 1$.
We first show that $r=1$ is a necessary condition for an infinite
number of solutions. 
Indeed if  $r \ge 2$, then we can eliminate variables $x_i = \pm 1$ and reduce
to an equation with smaller $n$  having all $|x_i| \ge 2$, with the same $(r,s)$ values
and a possibly
different value of $m$, call it $m'$. But now Theorem~\ref{th41} (1) rules out
$m'=0$, so we must have $m' \ne 0$.  For  $m' \ne 0$, 
Theorem~\ref{th51} gives an upper bound on the size of the solutions
which is independent of $m'$, thus establishing finiteness in this case.  

Suppose $r=1$.  Then the integer  solutions to the affine equation \eqn{121} 
consist of solutions to \eqn{601} plus additional integer solutions having
some variable $x_i=0$. We show the solutions with $x_i=0$ are finite
in number.  Indeed, if $x_n=0$ then \eqn{121} simplifies in the remaining variables to
$$
r x_1 x_2 \cdots x_{n-1} + s = 0.
$$
This clearly has finitely many solutions, since each $|x_i|\le |s|$. 
The theorem  now follows, using Theorem~\ref{th61} to classify all
cases when \eqn{601} has infinitely
many integer solutions. $~~~\bsq$

\end{document}